\documentclass[10pt,twosided]{article}
\usepackage[tbtags]{amsmath}
\usepackage{epsfig,amstext,amssymb,amsthm,latexsym}
\allowdisplaybreaks[4]
\usepackage{color}

\usepackage{amssymb,color}

\definecolor{c20}{rgb}{0.,0.7,0.}
\definecolor{c30}{rgb}{0.,0.,1.}
\definecolor{c40}{rgb}{1,0.1,0.7}
\definecolor{c50}{rgb}{1,0,0}
\definecolor{c60}{rgb}{1,0.9,0.1}

\def\yE#1{\textcolor{c30}{#1}}
\def\yE#1{#1}

\def\XE#1{\textcolor{c30}{#1}}
\def\XE#1{{#1}}

\def\enk#1{\textcolor{c30}{#1}}
\def\enk#1{#1}

\def\cE#1{#1}
\def\cT#1{#1}

\def\aH#1{\textcolor{c30}{#1}}
\def\aH#1{#1}

\def\dT#1{#1}
\def\dE#1{#1}

\def\eE#1{#1}

\def\cL#1{#1}

\def\cH#1{\textcolor{c20}{#1}}
\def\cH#1{#1}

\def\fD#1{\textcolor{c20}{#1}}
\def\fD#1{#1}

\def\gT#1{\textcolor{c40}{#1}}
\def\gT#1{#1}
\def\cl#1{\textcolor{c40}{#1}}
\def\cl#1{#1}

\def\LP#1{\textcolor{c30}{#1}}
\def\ZX#1{\textcolor{c20}{#1}}
\def\ee#1{\textcolor{c50}{#1}}
\def\ee#1{#1}
\def\LP#1{#1}
\def\ZX#1{#1}
\def\LPJ#1{\textcolor{c30}{#1}}
\def\LPJ#1{#1}
\def\LPJi#1{\textcolor{c30}{#1}}
\def\LPJi#1{#1}
\def\lp#1{\textcolor{c50}{#1}}
\def\lp#1{#1}
\def\clp#1{\textcolor{c50}{#1}}
\def\clp#1{#1}

\usepackage{amssymb,color}

\newcommand{\kb}[1]{\boldsymbol{#1}}
\newcommand{\vk}[1]{\kb{#1}}
\def\kal#1{{\cal{ #1}}}

\newcommand{\ve}{\varepsilon}

\newcommand{\abs}[1]{\lvert #1 \rvert}

\newcommand{\E}[1]{\mathbb{E}\left(#1\right)}

\newcommand{\CO}[1]{\mathbb{C}ov\left(#1\right)}

\newcommand{\pk}[1]{\mathbb{P} \left( #1 \right) }

\newcommand{\pb}[1]{\mbox{\rm$\vk{P}$}\Bigl \{#1 \Bigr \}}

\newcommand{\R}{\!I\!\!R}

\newcommand{\N}{\!I\!\!N}
\newcommand{\inr}{\in \R}

\newcommand{\inn}{\in \N}
\newcommand{\ldot}{,\ldots,}

\newcommand{\limit}[1]{\lim_{#1 \to   \infty}}

\newcommand{\BQN}{\begin{eqnarray}}
\newcommand{\EQN}{\end{eqnarray}}
\newcommand{\BQNY}{\begin{eqnarray*}}
\newcommand{\EQNY}{\end{eqnarray*}}

\newcommand{\BS}{\begin{sat}}
\newcommand{\ES}{\end{sat}}
\newcommand{\BT}{\begin{theo}}
\newcommand{\ET}{\end{theo}}
\newcommand{\BK}{\begin{korr}}
\newcommand{\EK}{\end{korr}}

\newcommand{\BD}{\begin{de}}
\newcommand{\ED}{\end{de}}

\newcommand{\BIT}{\begin{itemize}}
\newcommand{\EIT}{\end{itemize}}
\newcommand{\BDI}{\begin{description}}
\newcommand{\EDI}{\end{description}}

\newcommand{\BRM}{\begin{remarks}}
\newcommand{\ERM}{\end{remarks}}

\newcommand{\BEL}{\begin{lem}}
\newcommand{\EEL}{\end{lem}}

\newtheorem{theo}{Theorem}[section]
\newtheorem{sat}[theo]{Proposition}
\newtheorem{de}[theo]{Definition}
\newtheorem{lem}[theo]{Lemma}

\newtheorem{korr}[theo]{Corollary}
\newtheorem{remark}[theo]{Remark}
\newtheorem{remarks}[theo]{Remarks}

\newcommand{\nelem}[1]{{Lemma \ref{#1}}}

\newcommand{\netheo}[1]{{Theorem \ref{#1}}}
\newcommand{\nekorr}[1]{{Corollary \ref{#1}}}

\newcommand{\prooftheo}[1]{ \textsc{Proof of Theorem} \ref{#1} }

\newcommand{\proofkorr}[1]{\textsc{Proof of Corollary} \ref{#1}}

\newcommand{\COM}[1]{}
\newcommand{\US}[1]{\underset{#1}\sup}

\newcommand{\QED}{\hfill $\Box$}

\def\AL{\cE{\alpha}}

\def\BL{\cE{\beta}}

\def\ST{\sum_{i=1}^n \lambda_i X_i(t)}

\def\SIN{ \widetilde{\sigma}}

\def\tSIi{ \widetilde{\sigma_i}^2}

\def\tAi{ \cE{\widetilde{A_i}}}
\def\tAj{ \cE{\widetilde{A_j}}}

\def\tDi{ \cE{\widetilde{D_i}}}

\def\tAi{ \cE{{A_i}}}
\def\tAj{ \cE{{A_j}}}

\def\tDi{ \cE{{D_i}}}

\def\tN{\cT{\widetilde{N}}}

\def\TxT{T}
\def\tPSu{ \eE{\widetilde{\Psi}(u)}}
\def\th{\Lambda}
\def\tHu{ \cE{\th_\alpha,_\beta(u)}}

\def\tDHi{\cE{\widetilde{\sigma_{i}}^{2}}}

\def\gTt{\cE{g(t)}}
\def\gTT{\cE{g(T)}}
\def\cAL{\cE{\kal{C}_{\AL,\BL}}}
\def\cA{\cE{\kal{C}}}

\def\tMT{ \cE{\widetilde{M}}}
\def\XT{ \dE{\{X(t), t\in [0,T]\}}}
\def\XIT{ \dE{\{X_i(t), t\in [0,T]\}}}
\def\Fbm{ \aH{fBm} }
\def\sfbm{\aH{sub-fBm}}

\date{}

\def\njj{\cH{\mathbb{I}}}
\def\CCC{\aH{\kal{C}_0}}

\newcommand{\limitsupT}[1]{\underset{t\in[0,T]}\sup{#1}}
\def\longT{\mathcal{L}}
\def\heavyH{\mathcal{H}}
\def\IF{\infty}

\newcommand{\limitinf}[1]{\liminf_{#1 \to   \infty}}
\newcommand{\limitsup}[1]{\limsup_{#1 \to   \infty}}

\def\ZX#1{#1}

\def\LLP#1{\textcolor{c50}{#1}}
\def\LLP#1{#1}

\begin{document}
\title{\bf \Large  \dE{Finite-time Ruin \LLP{Probability} of \LLP{Aggregate}\\ Gaussian  Processes}}
\author{Krzysztof D\c{e}bicki\footnote{Mathematical Institute, University of Wroc\l aw, pl. Grunwaldzki 2/4, 50-384 Wroc\l aw, Poland},
Enkelejd Hashorva\footnote{Department of Actuarial Science, Faculty of Business and Economics,
University of Lausanne, Switzerland}, 
\cL{Lanpeng Ji,$^\dag$} 
 Zhongquan Tan\footnote{College of Mathematics, Physics and Information Engineering, Jiaxing University, Jiaxing 314001, PR China}
}
 \maketitle
\centerline{ \date{\today}}
 \baselineskip 15pt

\begin{quote}
{\bf Abstract:} Let  $\left\{\ST -\gTt, t\in [0,T]\right\}$ be an aggregate Gaussian risk process with a trend $\gTt$. We derive exact asymptotics of the finite-time ruin probability given by
$$\mathbb{P}\left(\sup_{t\in[0,T]}\left(\ST- \gTt  \right)>u\right)$$
as $u\to\IF$ for $\XIT, i\leq n,$ satisfying some asymptotic conditions. 
Further, we derive asymptotic results for the finite-time ruin probabilities of risk processes perturbed by an aggregate Gaussian process.

{\bf Key Words:} ruin probability; Gaussian process; perturbed risk process; \LP{L\'{e}vy process;}
\eE{(sub- and \gT{bi-})fractional Brownian motion;} risk aggregation; subexponential risks. 

{\bf AMS Classification:} Primary 60G15; Secondary 60G70, 68M20.

\end{quote}

\section{Introduction}
Numerous contributions have discussed the evaluation of the first-passage density of a random process $ \XT $ to a given deterministic
boundary denoted by $u+g(t)$ with fixed $u\ge0$. In a concrete insurance setup, \LP{let} $X(t)$ model the surplus process of the whole company \LP{at time $t$}, the decision to pay
dividends can be objectively made once the surplus process \ZX{crosses} the boundary. 
 Specifically, from the actuarial point of view,
it is of interest to calculate the crossing probability
\BQN\label{eq1111}
\pk{ \exists{t\in [0,T]},\ X(t)> u+g(t)}
\EQN
 for $u\ge 0$. However, an explicit formula for \eqref{eq1111} is hard to obtain except for some very special cases, e.g., $\{X(t), t\in[0,T]\}$ is a Brownian motion (Bm) and $g(t)$ is a linear function. Therefore, usually the aim of the analysis is to find adequate  approximations for it. \yE{From} risk theory point of view Eq. \eqref{eq1111} can also be seen as the finite-time ruin probability of an insurance company, i.e.,
\BQNY
\pk{ \exists{t\in [0,T]},\ X(t)> u+g(t)} =\ \pk{ \inf_{t\in [0,T]}\left(u+g(t)-X(t)\right)<0},
\EQNY
where $u\ge 0$ is the initial capital, $g(t)$ is the premium amount received up to time $t$, and $X(t)$ represents the aggregate claim amount up to $t$.
\LLP{Recently, the study of surplus process with dependent risks becomes more and more popular since independent risks is not applicable to practice, see e.g., Denuit et al. (2005) and  Constantinescu et al. (2011).}

In Michna (1998) \yE{it is shown that} the finite-time ruin probability given by
\BQN\label{eq1112}
\pk{\inf_{t\in [0,T]} \Bigl( u+ ct - B_H(t)\Bigr) < 0}
\label{rpfbm}
\EQN
\cH{is an \aH{adequate} approximation}
of the finite-time ruin probability for a  risk process with certain dependent risks, where $\{B_{H}(t), t\in [0,T]\}$ is a  fractional Brownian motion \aH{(\LLP{fBm})} with Hurst index $H\in (0,1]$.

\COM{
\ZX{For the asymptotic behaviour of supremum over
the infinite-time horizon, we mention the
seminal paper H\"{u}sler and Piterbarg (1999) which} derives  the exact asymptotic behaviour of
\dE{the infinite-time ruin probability}
\begin{eqnarray}
\label{eq1.2}
\mathbb{P}\left(\sup_{t\geq 0}\Bigl( X(t)-ct^{\theta} \Bigr)>u\right), \ \ \mbox{as}\ \ u\rightarrow\infty,
\end{eqnarray}
\cl{where} $\theta>H$ \LP{is a constant} and $\{X(t),t\in [0,\infty)\}$ \cl{is} a general Gaussian process.\\
}

Nowadays, all insurance  \cH{companies} run \cH{diverse} lines of business, \cH{with} typically some lines of business
(for non-life insurer) \lp{having} very high premiums because of high risks. \aH{In order to} reflect different portfolio variances,
as well as different business volumes, it is adequate to consider a  process which is a result of aggregation
of the specific portfolios. A tractable choice here is the \LLP{aggregate} process
\BQN\label{eqR}
 X(t) &= &\lambda_1 B_{H_1}(t)+ \cdots + \lambda_n B_{H_n}(t), \quad t\in [0,T],
 \EQN
where $\lambda_i,i\le n,$ are positive weights \aH{assigned} to \cH{the processes} $\{B_{H_i}(t), t\in [0,T]\}, \gT{i\le n},$ \XE{being} independent fBm's with Hurst indexes $H_i\in (0,1], i\le n$, respectively.\\
 \XE{Clearly, $X(t),t\in [0,T]$ is not a fBm} anymore; bounds and asymptotics of the finite-time ruin probability \XE{for $X(t),t\in [0,T]$} are given in D\c{e}bicki and Sikora (2011) for this multiplexed fBm's with a linear trend. \fD{ The asymptotics of the infinite-time ruin probability of the multiplexed fBm's with a trend is discussed in H\"{u}sler and Schmid (2006)}. 

\LLP{ The perturbed risk model 
is an important extension of the classical risk model.  
 Of course, instead of the Bm, general processes, including L\'{e}vy and Gaussian processes, can be considered as perturbations, see e.g., Schlegel (1998), Furrer (1998) and Frostig (2008).
In fact, the Bm (and L\'{e}vy processes) can not be
justified if the {\it perturbation terms} do not come from an i.i.d. framework, whereas some Gaussian processes can be.
In practice, the surplus is influenced by various uncertainties such as premium adjustments, legislation changes, cost of repairs, and other related expanses. Therefore, in order to reflect different variances of the uncertainties, it is reasonable to consider an aggregate Gaussian process as the perturbation.\\
}

\COM{
\LP{This} model is discussed \ZX{by} D\c{e}bicki and Sikora (2011).
 \cE{Specifically, for $\{\cH{B_{H_i}}(t), t\in [0,T]\}, i\le n,$ independent  \LLP{fBm}'s  with Hurst \cH{parameters} $H_i, i\le n,$ satisfying \dE{$0 < H:=H_1 < \cdots < H_n\le 1,$} the aforementioned  paper  derives the following exact asymptotic 
\BQN\label{eq:DS}
\pk{ \sup_{t\in [0,T]} \Bigl( \sum_{i=1}^n \lambda_i B_{H_i}(t) - c t\Bigr) > u} &\sim &  \CCC  \cH{u}
^{ \njj (H\le 1/2)\frac{1- 2 H}{H} }
\Psi\left(\frac{ u+ cT}{\eE{\sqrt{\sum_{i=1}^n \lambda_i T^{2 H_i}}}}\right),
\EQN
}
\dE{where $\Psi(u):=\pk{\eE{B_{1/2}(1)}>u},u\inr$, and $ \CCC $ is \ZX{a} known constant.}
Throughout this paper $\sim$ means the asymptotic equivalence of two functions when the argument tends to infinity, and $ \njj ( \cdot)$ stands for the indicator function.
}

\COM{
The results of Michna (1998) and D\c{e}bicki and Sikora (2011), which include upper and lower bounds for the finite-time ruin probability, are of interest for the dividend payment problem as discussed above.
Moreover, those results also have applications \ZX{to} the
pricing of certain barrier options, see e.g., Novikov et al.\ (2003). 
\\
}

\COM{
Let $\{\cE{X}(t),t\in[0,\infty)\}$ be a random processes with continuous sample paths. The problem of analyzing the asymptotics of
\begin{eqnarray}
\label{eq1.1}
\mathbb{P}\left(\sup_{t\in I}X(t)>u\right), \ \ \mbox{as}\ \ u\rightarrow\infty,
\end{eqnarray}
where $I$ is a finite or infinite interval, plays an important role in many fields of applied and theoretical probability.
For instance, the asymptotics behaviour of \eqref{eq1.1} is related to the asymptotics of ruin
probabilities in perturbed risk models in insurance and finance, see e.g., Embrechts et al.\ (1997).
For interesting applications related to Gaussian queueing models we refer the reader to D\c{e}bicki and Palmowski (1999),
Mandjes (2007),  D\c{e}bicki and Sikora (2011) and the references therein.  In fact, the study of the Gaussian model has been considered by many researchers. For $I=[0,\infty)$ and $\{X(t), t\in[0,\infty)\}$ a
Gaussian process with mean 0 and variance $t^{2H},0<H<1$, the seminal contribution  H\"{u}sler and Piterbarg (1999) derives  the exact asymptotic behaviour of the linear boundary crossing probability
\begin{eqnarray}
\label{eq1.2}
\mathbb{P}\left(\sup_{t\geq 0}(X(t)-ct^{\theta})>u\right), \ \ \mbox{as}\ \ u\rightarrow\infty
\end{eqnarray}
for any $\theta>H$.  Given the importance of the approximation of the ruin probabilities in finite time,
Michna (1998) and D\c{e}bicki and Rolski (2002) studied the finite-time ruin probability
\begin{eqnarray}
\label{eq1.3}
\mathbb{P}\left(\sup_{t\in [0,T]}(B_H(t)-ct)>u\right), \ \ \mbox{as}\ \ u\rightarrow\infty,
\end{eqnarray}
with $c>0, H\in (0,1],T>0$ and $B_H$ the  \LLP{fBm} with Hurst parameter $H$.
The finite-time ruin \cE{model} (\ref{eq1.3}) was \cE{generalised} by D\c{e}bicki and Sikora (2011) to the multiplex
case. \cE{Specifically, for $\{B_{H_i}(t), t\in [0,T]\}, H_i\in (0,1], i\le n$ independent fractional Brownian
 motions, the aforementioned  paper  derives the following exact asymptotics 
\BQN\label{eq:DS}
\pk{ \sup_{t\in [0,T]} \Bigl( \sum_{i=1}^n \lambda_i B_{H_i}(t) - c t\Bigr) > u} &\sim & {\cal C} \cH{u}
^{ \njj (H_1\le 1/2)\frac{1- 2 H_1}{H_1} }
\Psi\left(\frac{ u+ cT}{\sum_{i=1}^n \lambda_i T^{2 H_i}}\right),
\EQN
for any $c,\lambda_i,i\le n$ positive constants, where $\Psi(u)=\pk{B_{1}(t)>u},u\inr$, $\cal{C}$ is some known constant.
Throughout this paper $\sim$ means the asymptotic equivalence of two functions when the argument tends to infinity, and
$ \njj ( \cdot)$ stands for the indicator function. }
}

{In this paper we \aH{present} some extensions of D\c{e}bicki and Sikora (2011) \yE{and consider further the perturbed risk process}.
 \cH{Specifically,}  instead of dealing with the aggregation of independent  \LLP{fBm}'s, we consider the aggregation of independent centered Gaussian processes $ \XIT, i\le n$, \cl{with} some positive weights
$\lambda_i,i\le n $. Our analysis  then focusses on the asymptotics of \dE{the finite-time} ruin probability }
\BQNY
\pk{ \sup_{t\in [0,T]}  \Bigl(\sum_{i=1}^n \lambda_{i}X_i(t) - \gTt \Bigr) > u}, \ \text{\fD{as}}\ u\to\IF,
\EQNY
\eE{with some bounded measurable trend function \ZX{$\gTt$}.} It is worth noting that the aggregate Gaussian process $\sum_{i=1}^n \lambda_{i}X_i(t)$ is also a Gaussian process, but in order to see which of the components will contribute more to the asymptotics we would like to deal with the aggregate  Gaussian process other than one single Gaussian process. This might also be necessary from  practical point of view. \LLP{Moreover, the finite-time ruin probability of a perturbed risk process with perturbation modeled by an aggregate Gaussian process defined by
\BQNY
\pk{ \sup_{t\in [0,T]}  \Bigl(U(t)- c(t)+\sum_{i=1}^n \lambda_{i}X_i(t) \Bigr) > u}, \ \ {\fD{ u\ge0}},
\EQNY
is also discussed, where $U(t)-c(t)$ is the claim surplus process, and $\sum_{i=1}^n \lambda_{i}X_i(t)$ is the aggregate Gaussian perturbation. }

In the first result \netheo{thm1} we provide 
the asymptotic behaviour of the finite-time ruin probability \LLP{for the aggregate Gaussian process}, which indicates that the processes which have the smallest characteristic constants will contribute more to the asymptotics. 
Furthermore, \LLP{our second result \netheo{Th4.1} derives a novel asymptotic result for the finite-time ruin probabilities of some quite general perturbed risk processes including Gaussian perturbed risk process as a special case.}\\


\ZX{This paper is organized as follows.} In Section 2 we introduce some notation. 
The main results are given in Section \cl{3 and Section 4. Section 5 presents several examples}. \cH{Proofs} of all the results are relegated to Section 6.

\section{Notation and Preliminaries}
In this section we mention  several abbreviations and notation needed \cH{in this} paper \XE{and present the main assumptions}.
There are mainly two well known constants, \XE{namely} Pickands constant and  Piterbarg constant, which play important roles in the extreme theory  of Gaussian processes. The former is defined by
$$\mathcal{H}_{\alpha/2}=\lim_{T\rightarrow\infty} T^{-1} \E{ \exp\biggl(\sup_{t\in[0,T]}\Bigl(\sqrt{2}B_{\alpha/2}(t)-t^{\alpha}\Bigr)\biggr)},\quad \eE{\alpha\in (0,2],}$$
and the latter is  defined by
$$\mathcal{P}_{\alpha}^{R}:=\cE{\lim_{S \to \infty}}\E{ \exp\biggl(\sup_{t\in[0,S]}\Bigl(\lp{\sqrt{2}}B_{\alpha/2}(t)-(1+R)t^{\alpha}\Bigr )\biggr)},\ \ \alpha\in(0,2], \ R>0,$$
where $\{B_{\alpha/2}(t),t\in[0,\IF)\}$ is a fBm with Hurst index $\alpha/2$.
See  \enk{Pickands (1969)} or  Piterbarg (1996), 
for the main properties of Pickands and Piterbarg constants.\\

\yE{We shall impose two main common assumptions on the Gaussian processes of interest.} \cH{
\LLP{Let $\{\xi(t),t\in [0,\IF)\}$ be a centered Gaussian process with variance function $\sigma_\xi^2(\cdot)$. Throughout this paper \fD{the} process $\xi$ with a bar represents a standardized process i.e., $\bar{\xi}(t):= \xi(t)/\sigma_\xi(t)$}}.


\cE{{\bf Assumption A1}.  The standard deviation function $\sigma_\xi(\cdot)$ of the Gaussian process $ \xi(t) $ attains
its maximum, denoted by $\widetilde{\sigma}$, over $[0,T]$ at the unique point $t=T$.  Further, \fD{there exist} some positive constants
$\fD{\alpha\in (0,2],\beta},A,D$ such that
\BQN
\eE{\sigma_\xi(t)}&=& \SIN -A(T-t)^{\beta}+o((T-t)^{\beta}), \quad t\to T,
\EQN
}
and
$$\eE{\CO{\bar{\xi}(s), \bar{\xi}(t)}=1- D|t-s|^{\alpha}+o(|t-s|^{\alpha})}, \quad \min(t,s)\rightarrow T.$$

\cE{
{\bf Assumption A2}. There exist \ZX{positive constants $\mathbb{C}, \delta$ and $\gamma$} such that, for all $s,t\in[\delta, T]$,
\BQN
\E{(\xi(t)-\xi(s))^{2}} &\leq & \mathbb{C}|t-s|^{\gamma}.
\EQN
}

Some recent \cH{studies} in financial markets indicate that the class
of $H$-self-similar \aH{($H$-ss)} Gaussian \aH{processes}
can \aH{adequately} model the long-range dependence structure of the real financial data.
 Let us recall that
 a centered Gaussian process $\{X(t),t\in[0,\infty)\}$ with $X(0)=0$ is  $H$-ss with an exponent $H\in (0,1\cH{]}$ if the covariance \eE{function satisfies the condition}
$$\CO{X(at),X(as)}=a^{2H}\CO{X(t),X(s)}, \quad \aH{\forall a\in (0,\infty)}.$$
\COM{
\aH{Furthermore,} if \LP{the $H$-ss Gaussian \aH{process}} $\{X(t),t\in[0,\infty)\}$ is a martingale, then
$$\CO{\eE{\LP{X}(t),\LP{X}(s)}}=\LP{\E{X^2(1)}}\min\aH{(}s^{2H},t^{2H}),$$
which further yields that
$$
X(t)=\sqrt{\E{X^2(1)}}B_{1/2}(t^{2H}),
$$
a time-space-changed Bm. Additionally,
}
\eE{\XE{A} \aH{prominent example of} self-similar Gaussian processes is the \cH{bi-}fractional  Brownian motion \aH{(bi-\LLP{fBm})} $\{B_{K,H}(t),t\in[0,\infty)\}$ \aH{with covariance function given by}}
\BQNY
\CO{B_{K,H}(t),B_{K,H}(s)}=\frac{1}{2^{K}}[(t^{2H}+s^{2H})^{K}-|s-t|^{2KH}],\quad K\in(0,1], \quad \LP{H\in(0,1).}
\EQNY
\eE{Another interesting self-similar Gaussian process is the sub-fractional Brownian motion \aH{(sub-fBm)} $\{S_H(t),t\in[0,\infty)\}$ with covariance function given by
\BQNY
\CO{S_H(t),S_H(s)}=t^{2H}+s^{2H}-\frac{1}{2}\left[(s+t)^{2H}+|t-s|^{2H}\right], \quad \LP{H\in(0,1).}
\EQNY
\cl{
\cH{\aH{Important} results} \aH{for}  \cH{the} \cH{bi-}\LLP{fBm} \aH{and sub-\LLP{fBm}} can be found in Houdr\'{e} and Villa (2003) \ZX{and }
Bojdecki et al.\ (2004).}
}

\COM{

We conclude this section with Piterbarg-Prisyazhnyuk Theorem (see Theorem 1 of \lp{Piterbarg and Prisyazhnyuk} (1978)).

\BT \label{thm:pit} Let $\XT$ be \cE{a} centered  Gaussian process with \cH{a.s.}\ continuous  sample paths  satisfying Assumption A0 with $\eE{\sigma_X(T)}=1$. Assume that
\ee{Assumption A1 holds with $\widetilde{\alpha},\dT{\widetilde{\beta}},A,D$, and further Assumption A2 is satisfied
for some positive constants $C,\gamma$ and $\eE{\delta}$.}

 Then, \ZX{as $u\to \infty$, we have}
\BQN
\mathbb{P}\left(\sup_{t\in[\LPJi{\delta},T]}X(t)>u\right)&\sim & \LPJi{\mathcal{C}}_{\dT{\widetilde{\alpha}},\dT{\widetilde{\beta}}}     \th_{\dT{\widetilde{\alpha}}, \dT{\widetilde{\beta}}}(u) \Psi(u),
\EQN
where
\BQNY
\cA_{\dT{\widetilde{\alpha}}, \dT{\widetilde{\beta}}}=\left\{
              \begin{array}{ll}
\mathcal{H}_{\dT{\widetilde{\alpha}}/2}\Gamma(1/\dT{\widetilde{\beta}})D^{1/\dT{\widetilde{\alpha}}}\dT{\widetilde{\beta}}^{-1}A^{-1/\dT{\widetilde{\beta}}}     , & \hbox{if } \dT{\widetilde{\alpha}}< \dT{\widetilde{\beta}} ,\\
\mathcal{P}_{\dT{\widetilde{\alpha}}}^{A/D}, & \hbox{if } \dT{\widetilde{\alpha}} =\dT{\widetilde{\beta}},\\
1, & \hbox{if } \dT{\widetilde{\alpha}} >\dT{\widetilde{\beta}},
              \end{array}
            \right.
\quad\mathrm{and}\quad
\th_{\dT{\widetilde{\alpha}},\dT{\widetilde{\beta}}}(u)=\left\{
              \begin{array}{ll}
u^{2/\dT{\widetilde{\alpha}}-2/\dT{\widetilde{\beta}}},& \hbox{if } \dT{\widetilde{\alpha}} <\dT{\widetilde{\beta}},\\
1, & \hbox{if } \dT{\widetilde{\alpha}} \ge\dT{\widetilde{\beta}}.
              \end{array}
            \right.
\EQNY
\COM{
(i) for $\beta>\alpha$ and $\mathcal{G}_{\alpha,\beta}=\mathcal{H}_{\alpha/2}\Gamma(1/\beta)D^{1/\alpha}\beta^{-1}A^{-1/\beta}$, as $u\rightarrow\infty,$
$$\mathbb{P}\left(\sup_{t\in[\eE{\delta},T]}\xi(t)>u\right)=\mathcal{G}_{\alpha,\beta}u^{2/\alpha-2/\beta}\Psi(u)(1+o(1));$$

(ii) for $\beta=\alpha$ and $R:=A/D$, as $u\rightarrow\infty,$
$$\mathbb{P}\left(\sup_{t\in[\eE{\delta},T]}\xi(t)>u\right)=\mathcal{P}_{\alpha}^{R}\Psi(u)(1+o(1));$$

(iii) for $\beta<\alpha$, as $u\rightarrow\infty,$
$$\mathbb{P}\left(\sup_{t\in[\eE{\delta},T]}\xi(t)>u\right)=\Psi(u)(1+o(1)).$$
}
\ET

}

\section{\LLP{Exact Asymptotics of the Finite-time Ruin Probability}}

Given $n$ independent centered Gaussian processes $ \XIT, i\le n,$ with a.s. continuous sample paths and standard deviation \fD{functions} $\sigma_i(\cdot), i\le n$, respectively,
the extended \dE{D\c{e}bicki-Sikora} Gaussian model consist\eE{s} in the specification of the \LLP{aggregate} Gaussian process
\BQN\label{mod}
 X(t):= \lambda_1 X_1(t) + \cdots + \lambda_n X_n(t), \quad  t\in [0,T],
\EQN
with \LP{$\lambda_i\ge0, i\le n$}. The finite-time ruin probability of this risk model is defined as
\BQNY
\pk{ \sup_{t\in [0,T]}  \Bigl(X(t) - \gTt \Bigr) > u},
\EQNY
 for the deterministic \LLP{bounded measurable} trend function $\gTt$ and $u\ge 0$.\\  

\COM{
\cE{\aH{For the dividend problem discussed in Teunen and Goovaerts (1993), as well as other interesting problems related to pricing of barrier options it is crucial} 
to derive sharp lower and upper \LLP{ bounds} for all values of the parameter (in our case $u$). \aH{Our first result 
 provides useful lower and upper bounds for the \XE{boundary} crossing probability of interest}.}
In the following, \eE{assume that $h^*(t):=\sigma_X^2(t)$ is a strictly increasing function on $[0,T]$, and denote \ZX{its inverse function by $h$}. }
\gT{For a given trend function $g$, we denote the polygonal line with $k$ nodes $(t_0,g(h(t_0))) \ldot (t_k,g(h(t_k)))$ by $\overline{g_k}$, 
where $0=t_0 < t_1 < \cdots < t_k=\eE{h^*(T)}$.} \LP{Furthermore, we} denote by $\eE{\Gamma_T}$
the set of all such polygonal lines with nodes defined by $g$ \LP{and $h$}.

\BT \label{thm2} Let $ \XIT, i\leq n,$ be independent centered Gaussian processes with a.s. continuous sample paths. 
 \cH{We have}, for $g\aH{(t)}$ some bounded measurable trend function and $u \LP{\ge}0$,
\BQN
\mathbb{P}\left(\sup_{t\in[0,T]}\left(\sum_{i=1}^{n}\lambda_{i}X_{i}(t)- g(t)  \right)>u\right)
& \geq & \dE{\sup_{t\in [0,T]}}\cT{\Psi\left(\frac{u+ g(t)  }{\sqrt{\eE{h^*(t)} } }\right).}
\EQN
 If each \LP{of} $\XIT, i\leq n,$ is \cH{either} \eE{a \Fbm with parameter $H_i\in[1/2,\gT{1]}$, \cH{or} a \sfbm \   with parameter \LLP{$H_i\in(0,1/2],$} or a \LLP{Gaussian martingale with increasing variance function}, then we have 
 }
%
\ZX{\BQN\label{eq:last}
\pk{\sup_{t\in[0,T]}\left(\sum_{i=1}^{n}\lambda_{i}X_{i}(t)- g(t)  \right)>u}
\le 1-\eE{N_T(u,g,h)},
\EQN}
provided that $g(h(t))$ is a concave function \gT{on} $[0, \eE{h^*(T)}]$, where
\BQNY
N_T(u,g,h)&=&\sup_{\overline{g_k} \in \Gamma_T}
\aH{\Biggl( }
\frac{1}{\sqrt{(2\pi)^{k} \prod_{i=0}^{k-1}(t_{i+1}-t_i) }}
\int_{-\infty}^{u+ \gT{g(h(t_k))}} \cdots \int_{-\infty}^{u+ \gT{g(h(t_1))}}
\exp\left(-\sum_{i=1}^ {k-1} \frac{ (x_{i+1}- x_i)^2}{2(t_{i+1}- t_i)} \right)\notag \\
&& \qquad\times\prod_{i=0}^{k-1}
\left[1-\exp\left(-2\frac{ ( u+g(h(t_i)) - x_i)(u+g(h(t_{i+1}))-x_{i+1})}{ t_{i+1}-t_i} \right)\right]\,  dx_1 \cdots d x_k
\aH{\Biggr)}.
\EQNY
\ET

\BRM \label{remarks:12}
(a) \LP{If each of $\XIT, i\leq n,$ is \eE{a \LLP{fBm}  or a  $H$-ss} Gaussian martingale, with $H_i\in[1/2,1], i \le n$, then the inverse function $h$ is concave. Moreover, if the trend function $g$ is concave and increasing, the concavity condition on $g(h(t))$ in \netheo{thm2} is fulfilled. }

(b) Typically, \ZX{in insurance applications it is} assumed that premium rate $c>0$ is constant. Consider the processes in (a).
In view of \eqref{eq:last},  
\ZX{for any positive $c, T$ and $u$ we have}
\BQN\label{eqb314}
{\pk{\sup_{t\in[0,T]}\left(\sum_{i=1}^{n}\lambda_{i}X_{i}(t)- ct  \right)>u}}\le 1-{N_T(u,g,h)}\le \cT{ \Psi\left(\frac{u+ c T   }{\sigma_X(T)}\right)
+\exp\left(\frac{-2 c T u}{\sigma_X^2(T)} \right)\Psi\left(\frac{u- c T   }{\sigma_X(T) }\right)}.
\EQN
\aH{Note that \ZX{ the right-hand
side \ZX{of \eqref{eqb314}} \LP{has been} given in Theorem 2.4 of D\c{e}bicki and Sikora (2011) for the case of \LLP{fBm}'s.}}
\ERM

}

In order to obtain the exact asymptotics of the finite-time ruin probability, some conditions on the Gaussian processes and the \fD{bounded measurable trend function $\gTt$} \LLP{needed are} fully described in \netheo{thm1}. \dE{
\cl{\aH{For our results below we need the following notation}
\BQNY
\tHu:=\left\{
              \begin{array}{ll}
\left(  \frac{u+ g(T)}{\eE{\sqrt{\sum_{i=1}^n \lambda_i^2 \tSIi}}}\right)^{\cT{2/\AL-2/\beta}}, & \hbox{if } \dT{{\alpha}} < \dT{{\beta}},\\
1, & \hbox{if } \dT{{\alpha}} \ge\dT{{\beta}},
              \end{array}
            \right.
\quad\mathrm{with}\quad\widetilde{\sigma}_i:=\sigma_i(T).
\EQNY
}}
\LLP{Further,  $\Gamma(\cdot)$ stands for the Euler Gamma function and $\njj(\cdot)$ for the indicator function.} 
\LPJi{Next} we state our \LLP{first} result.
\BT \label{thm1} Let $ \XIT, i\le n,$ be independent centered Gaussian processes with a.s. continuous sample paths and standard deviation functions \gT{$\sigma_{i}(\cdot), i\le n$,} and define $ \XT $ as in \eqref{mod}.
If  Assumptions A1 and A2 hold \ZX{for each $\XIT, i\le n,$} with constants
 $ \alpha_i,\beta_i,  A_i,D_i,\mathbb{{C}},\delta, \gamma_i, i\le n$, respectively,   \clp{then, for any  bounded measurable trend function $\gTt$ satisfying
\BQN\label{eqgg}
\big|g(T)-g(t)\big|\le\mathcal{M}( T-t)^{\min_{ i\le n}\beta_i},\ \  \forall t\in[\nu, T] 
\EQN
 for some constant $\mathcal{M}$ and $\nu\in(0,T)$,} 
  we have
\BQN \label{eq:thm:1}
\mathbb{P}\left(\sup_{t\in[0,T]}\left( \LLP{X(t)}- \gTt  \right)>u\right)
&\sim&  \cAL      \tHu  \Psi\left( \frac{u+ g(T)}{\cL{\sqrt{\sum_{i=1}^n \lambda_i^2 \tSIi}}} \right), \quad u\rightarrow\infty,
\EQN
where
$$ \cAL    =\left\{
              \begin{array}{ll}
                \mathcal{H}_{{\AL}/{2}}\LLP{\Gamma(1/\beta+1)}\LLP{\tN^{-1/\beta}\widetilde{G}^{1/\AL}\left({\sum_{i=1}^{n}\lambda_{i}^{2} \tDHi}\right)^{1/\beta-1/\alpha} }, & \hbox{if } \AL <\cT{\beta}, \\
\eE{\mathcal{P}_{\alpha}^{\widetilde{N}/\widetilde{G}} }
                , & \hbox{if } \AL =\cT{\beta},\\

1 , & \gT{\hbox{if } \AL >\beta, }

              \end{array}
            \right.
$$
\lp{with}
$$\LLP{\alpha=\min_{i\le n}\alpha_i,}\  \beta=\min_{ i\le n}\beta_i, \ \cT{\widetilde{N}=\sum_{i=1}^{n} \lambda_{i}^{2}\widetilde{\sigma_{i}}\LLP{A_i\njj(\beta_{i}=\beta)}},\ \widetilde{G}=\sum_{i=1}^{n}\lambda_{i}^{2}\tDi \LLP{\widetilde{\sigma_{i}}^2\njj(\alpha_{i}=\AL )}.$$
\ET

\BK \label{korr:1}
\yE{Let $\XIT, i\leq n$, $\XT$ and \fD{ $\gTt$} be as in \netheo{thm1}. \\
(i) If $\XIT, i\leq n,$ are \eE{\aH{bi-\LLP{fBm}'s} with parameters $K_i,H_i\in(0,\gT{1]}$, $i\le n,$ satisfying $0< \dE{KH}:=K_1H_{1}<K_2H_{2}\le\cdots\le K_nH_{n}$,} 
then we have}
\BQN
\mathbb{P}\left(\sup_{t\in[0,T]}\left(\LLP{X}(t)- \gTt  \right)>u\right)
&\sim &
\gT{\cA_{2 KH, 1}\th_{2 KH, 1}(u)\Psi\left( \frac{u+ g(T)}{\cL{\sqrt{\sum_{i=1}^n \lambda_i^2 \tSIi}}} \right)},  \quad u\to \infty,
\EQN
where
$$\eE{\cA_{2 \eE{KH}, 1}=\left\{
            \begin{array}{ll}
\mathcal{H}_{\dE{KH}}\left(\sum_{i=1}^{n} \lambda_{i}^{2}\LP{\widetilde{\sigma_{i}}^2}\right)^{\frac{2KH-1}{2\aH{KH}}}\frac{(\frac{1}{2^{K}}\lambda_{1}^{2})^{1/(2\aH{KH})} \TxT}{\sum_{i=1}^{n} \lambda_{i}^{2}K_i\cE{H_{i}}\LP{\widetilde{\sigma_{i}}^2}}            , & \hbox{if } \dE{\aH{KH}} <1/2 ,\\
1+\frac{\lambda_{1}^{2}\TxT}{\lp{{2^{K}}(\sum_{i=2}^{n}\lambda_{i}^{2}K_i H_i \LP{\widetilde{\sigma_{i}}^2}+\lambda_{1}^{2}\TxT/2)}}, & \hbox{if } \dE{KH} =1/2,\\
1, &\hbox{if } \dE{KH}>1/2,
              \end{array}
            \right.
\ \ \mbox{and}\ \  \gT{\widetilde{\sigma_{i}}}=T^{K_iH_{i}}.}$$

(ii) \eE{If $\XIT, i\leq n,$ are sub-\LLP{fBm}'s with parameters $H_{i}\in (0,1), i\le n,$ satisfying $H:=H_{1}<H_{2}\le\cdots\le H_{n}$,
then
\BQN
\mathbb{P}\left(\sup_{t\in[0,T]}\left(\LLP{X}(t)- \gTt  \right)>u\right)
&\sim & \cH{ \cA_{2 H, 1}
\th_{2 H, 1}(u)\Psi\left( \frac{u+ g(T)}{\cL{\sqrt{\sum_{i=1}^n \lambda_i^2 \tSIi}}} \right)}, \quad u\to \infty,
\EQN
where
$$\cA_{2 \dE{H}, 1}=\left\{
            \begin{array}{ll}
\mathcal{H}_{\dE{H}}\left(\sum_{i=1}^{n} \lambda_{i}^{2}\aH{\widetilde{\sigma_{i}}^{2}}
\right)^{\frac{2H-1}{2H}}\frac{(\frac{1}{2}\lambda_{1}^{2})^{1/\cE{(2\dE{H})}} \TxT}{\sum_{i=1}^{n}\lambda_{i}^{2} H_i \aH{\widetilde{\sigma_{i}}^{2}} 
}            , & \hbox{if } \dE{H} <1/2 ,\\
1+\frac{\lambda_{1}^{2}\TxT}{ \lp{2\sum_{i=2}^{n}}\lambda_{i}^{2}H_i \aH{\widetilde{\sigma_{i}}^{2}} 
+\lambda_{1}^{2}\TxT}, & \hbox{if } \dE{H} =1/2,\\
1&  \hbox{if } \dE{H} >1/2,
              \end{array}
            \right.\ \ \mbox{and}\ \ \widetilde{\sigma_{i}}^{2}=\LP{(2-2^{2H_{i}-1})}T^{2H_{i}}.
$$
}
\COM{
(iii) If $\XIT, i\leq n,$ are  $H$-ss Gaussian \gT{martingales} with 
parameters $H_i\in(0,\gT{1]}, i\le n$,
then, 
as $u\to \infty$,
\BQN
\dT{\mathbb{P}\left(\sup_{t\in[0,T]}\left(\sum_{i=1}^{n}\lambda_{i}X_{i}(t)- \gTt  \right)>u\right)}
&\sim &
\cH{\lp{2}\gT{\Psi\left( \frac{u+ g(T)}{\cL{\sqrt{\sum_{i=1}^n \lambda_i^2 \tSIi}}} \right)}},
\EQN
where $\gT{\widetilde{\sigma_{i}}=\LP{\sqrt{\E{X_i^2(1)}}}T^{H_{i}}}.$}

\EK




\section{\LP{ Perturbed} Risk Processes}
This section 
is devoted to
the analysis of \cH{finite-time} ruin probabilities of some \LP{general perturbed risk models.}
In particular, we focus on perturbed risk processes, where the perturbation is an \LLP{aggregate} centered Gaussian \gT{process} \LLP{representing the aggregation of different types of perturbations}.
\LP{Consider the claim surplus process of an insurance company} defined by
\BQN
S(t)&=&\LLP{U}(t)-c(t), \quad t\ge0,\label{rk1}
\EQN
where $\{\LLP{U}(t), t\in[0,\infty)\}$ \COM{$(t)=\sum_{i=1}^{N(t)}Z_i$} is \LP{the}  \LLP{aggregate} claim process and $c(t)$  is a nonnegative increasing function \aH{modeling} the premium income. 
Further, define the claim surplus process of the perturbed risk process as
\BQN
\tilde{S}(t)&=&S(t)+X(t), \quad t\ge0, \label{prk2}
\EQN
where \lp{the process} $\{X(t), t\in[0,\infty)\}$\COM{= \sum_{i=1}^n \lambda_i X_i(t), t\ge0, is the \LLP{aggregate} centered Gaussian process} 
is a perturbation \LP{assumed to be} independent of $\{S(t), t\in[0,\infty)\}$.

For any \gT{$T\in(0,\infty)$,} \LP{the \cH{finite-time} ruin probability for the processes (\ref{rk1}) and (\ref{prk2}) are defined as}
\BQNY
\psi(u,c,T)=\pk {\underset{t\in[0,T]}\sup S(t)>u}\quad\mathrm{and}\quad\tilde{\psi}(u,c,T)=\pk {\underset{t\in[0,T]}\sup \tilde{S}(t)>u},
\EQNY
respectively, \gT{where $u\ge0$ is the initial surplus.}
\gT{In general, the calculation} of the \cH{finite-time} ruin probability is more difficult than the \cH{infinite-time} ruin probability.
\cH{Therefore, often the aim of the analysis is to find good approximation for it.}  For notational simplicity set below 
\BQNY
F_1(u)=\pk{\limitsupT{\LLP{U}(t)}\le u}, \quad F_2(u)=\pk{\limitsupT{X(t)}\le u}, \quad u\ge 0.
\EQNY

\COM{\cH{We consider} the Gaussian perturbed Sparre Andersen risk model where $c(t)=ct$, $\tau_i, i\inn$ are identically distributed random variables with common mean $\frac{1}{\delta}$, and $Z_i,i\inn,$ are independent identically distributed random variables.  \cH{Typically, in insurance applications net profit condition}
\BQN\label{profit}
c&>&\delta \E {Z_1}
\EQN
is assumed.
}
Let us first recall the class of long-tailed distributions and that of heavy-tailed distributions.

{\it Heavy-tailed distribution class $(\heavyH)$}: A distribution function $F$ is said to be heavy-tailed if and only if
\BQNY
\int_{\LLP{-\IF}}^\infty e^{\lambda x}F(dx)=\infty \quad \mathrm{for}\ \mathrm{all}\ \lambda>0.
\EQNY
{\it Long-tailed distribution class $(\longT)$}: A distribution function $F$ 
is said to be long-tailed \LLP{if and only if  
\BQNY
\limit{x}\frac{1-{F}(x+y)}{1-{F}(x)}=1 \quad \mathrm{for}\ \mathrm{all}\ y\in \R.
\EQNY
}

\COM{
\cH{{In} our asymptotic analysis, we shall consider the claims \aH{from} the class $\mathcal{S}$  of subexponential distributions.}
By definition, a distribution function $F$ supported on $[0,\infty)$ belongs to $\mathcal{S}$ if
\BQNY
\underset{x\rightarrow \infty}\lim\frac{\aH{1- }F^{*2}(x)}{\aH{1- }F(x)}&=&2,
\EQNY
where \aH{$F^{*2}$ denotes the $2$}-fold convolution of $F$; }

It is \ee{well-}known that $ \longT\subset\heavyH$, \ee{see e.g.,} Embrechts et al.\ (1997) and Foss et al.\ (2011) for the basic
properties of heavy-tailed distributions. In addition, $F\in\longT$ implies that there exists some function $d(u), u\ge 0$ such that
$$\lim_{u\to \IF} \frac{u}{d(u)}= \lim_{u\to \IF} {d(u)}= \IF$$
 and
 \BQN
1-{F}(u+ d(u))\ \sim\  1-{F}(u), \quad u\rightarrow\IF,\label{eq_2}
\EQN
see e.g., 
 Foss et al.\ (2011). Next, we present the main result of this \yE{section.}\\

\BT \label{Th4.1} 
\LLP{Assume that}  $F_1\in\longT$ and
 $1-{F_2}(u)=o(1-{F_1}(u))$ as $u\rightarrow \infty$, then
 \BQN
\tilde{\psi}(u,c,T)\ \sim\ 1-F_1(u)\ \sim \  {\psi}(u,c,T), \quad u\rightarrow\IF.\label{eqTh4.1_2}
\EQN
\ET

\LLP{In the following,} we consider Gaussian perturbed L\'{e}vy risk processes, where the perturbation is an \LLP{aggregate} Gaussian process  $X(t)= \sum_{i=1}^n \lambda_i X_i(t), t\ge0$,  discussed in Section 3. 

\COM{

 see for example, Berman (1986), Marcus (1987) Wilekens (1987), Albin (2007) and  Albin and Sund\'{e}n (2009).

\textbf{Proposition A.} Let $\{Y(t), t\ge0\}$ be a L\'{e}vy process with characteristic triple $(d, \sigma^2, \Pi)$. We have
\BQNY
\frac{\Pi([1,\IF)\cap\cdot)}{\Pi[1,\IF)}\in \longT \Rightarrow Y(T)\in \longT\ \mathrm{for}\ T>0\Leftrightarrow \limitsupT Y(t)\in \longT.
\EQNY
Moreover,
\BQNY
\pk{\limitsupT Y(t)>u}\ \sim\ \pk{Y(T)>u}\ \mathrm{as}\ u\rightarrow \IF.\frac{\Pi([1,\IF)\cap\cdot)}{\Pi[1,\IF)}\in \longT,
\EQNY
}
\BK\label{k4.2}
 If $\{\LLP{U}(t), t\in[0,\IF)\}$ is a L\'{e}vy process 
 such that
\BQN
\LLP{U}(T)\in \longT,\label{eqpi}
\EQN
and $\{X(t), t\in[0,T]\}$ is an \LLP{aggregate} Gaussian process satisfying the conditions of Theorem 3.1, then
\BQNY
\tilde{\psi}(u,c,T)\ \sim\ \pk{\LLP{U}(T)>u}\ \mathrm{as}\ u\rightarrow \IF.
\EQNY
\EK
\begin{remark}
\ee{In the light of Albin and Sund\'{e}n (2009)}, for a L\'{e}vy process $\{Y(t), t\in[0,\IF)\}$ with characteristic triple $(d, \sigma^2, \Pi)$,
\BQNY
\frac{\Pi([1,\IF)\cap\cdot)}{\LPJi{\Pi([1,\IF))}}\in \longT
\EQNY
\ee{implies that $Y(T)\in \longT$.}
\end{remark}
\COM{
\subsection{\LP{Bounds of The Finite-time Ruin Probability}}
\LP{This subsection discusses a Gaussian perturbed risk model with dependent claims. Consider a claim counting process $N(t)=\sum_{i=1}^{\infty}\mathbb{I}(\tau_1+\cdots+\tau_i\le t)$, with $\tau_i, i\inn,$ being \gT{independent} claim inter-arrival times and satisfying $\pk{\tau_i>y}=e^{-\delta_i y}$, $\delta_i>0$. \aH{Let}
 $Z_i,  i\inn,$ be the non-negative claim sizes which are integer-valued and assumed to be independent of
$\{N(t), t\ge0\}$. In such a case, $\LLP{U}(t)=\sum_{i=1}^{N(t)}Z_i$ is the \LLP{aggregate} claim process. Moreover, let $X(t)= \sum_{i=1}^n \lambda_i X_i(t), t\ge0$.}
Set $p_{z,i}=\pk{Z_1=z_1,\cdots,Z_i=z_i}$ with $z_i\inn$ denoting the joint distribution of $Z_j, j\le i$. Let $n_u^c=[u+c(T)]+1,$ $d_u^{c}(y)=\inf\{t;u+c(\gT{t})\ge y\}$, and $v_{u,i}^c=d_u^{c}(i), i\in \{0\}\cup \N$, where $[c]$ denotes the integer part of $c$. Let $k_{u}^c=k(z_1,\cdots,z_{n_u^c})$,
as a function of $z_1,\cdots,z_{n_u^c}$, \cH{be \LP{an} integer} such that
\BQNY
z_1+\cdots+z_{k_u^c-1}\le n_u^c-1, \quad \mathrm{and}\quad z_1+\cdots+z_{k_u^c}\ge n_u^c.
\EQNY
Furthermore, denote by $\mathcal{G}_T$ the set of functions $g$ such that $g(h(t))$ is concave
on $[0,h^*(T)]$ and $c(t)-g(t)$ is a positive increasing function
with $\underset{t\rightarrow\infty}\lim(c(t)-g(t))=\infty$. By $\mathfrak{G}_T$
we denote the set of functions $f$ such that $-f(h(t))$ is concave on $[0,h^*(T)]$ and $c(t)-f(t)$ is a positive increasing function with $\underset{t\rightarrow\infty}\lim(c(t)-f(t))=\infty$.

\gT{In Theorem \ref{Th5.4} we obtain  bounds for the \cH{finite-time} ruin probability which are
valid for any value of the initial surplus}.\\
\BT \label{Th5.4}
If each \LP{of} $\XIT, i\le n,$ satisfies conditions of \netheo{thm2}, then
\BQN
\US{f\in\mathfrak{G}_T}\US{\eta\in(0,\infty)}\left(\Bigl(1-L_T(u+\eta,c-f,\rho)\Bigr)N_T(\eta,-f,h)\right)\le\tilde{\psi}(u,c,T)\le 1-\US{g\in\mathcal{G}_T}\US{\xi\in(0,u)}\biggl(L_T(u-\xi,c-g,\mu)N_T(\xi,g,h)\biggr), \quad
\EQN
where
\begin{align*}
L_T(u,c,s)=\sum_{n_u^c=1}^\infty \gT{\mathbb{I}}(v_{u,n_u^c-1}^c\le T\le v^c_{u,n_u^c})\sum_{z_i\ge 1, i\le n_u^c}p_{z,n_u^c}e^{-s_{k_u^c} T}\frac{\prod_{i=1}^{k_u^c} \delta_i}{s_{k_u^c}^{k_u^c}}
\sum_{j=0}^{k_u^c-1}(-1)^j b_j(w^c_{u,1},\cdots,w^c_{u,j})s_{k_u^c}^j\sum_{m=0}^{k_u^c-j-1}\frac{(s_{k_u^c} T)^m}{m!},
\end{align*}
with $\rho_{k_u^c}=\min(\delta_1,\cdots,\delta_{k_u^c})$, $\mu_{k_u^c}=\max(\delta_1,\cdots,\delta_{k_u^c})$, $w^c_{u,l}=v^c_{u,z_1+\cdots+z_l}, l\inn$, and
\BQNY
b_0=1,\quad b_1(x_1)=x_1,\quad
b_j(x_{1},\cdots,x_{j})=\det\left(
                              \begin{array}{cccccc}
                                \frac{x_1}{1!} & 1 & 0 & 0 & \cdots & 0 \\
                                \frac{x_2^2}{2!} & \frac{x_2}{1!} & 1 & 0 & \cdots & 0 \\
                                \frac{x_3^3}{3!} &  \frac{x_3^2}{2!} &  \frac{x_3}{1!} & 1 & \cdots & 0 \\
                                \vdots & \vdots & \vdots & \vdots & \cdots & \vdots \\
                                  \frac{x_{j-1}^{j-1}}{(j-1)!}&  \frac{x_{j-1}^{j-2}}{(j-2)!} &  \frac{x_{j-1}^{j-3}}{(j-3)!} &  \frac{x_{j-1}^{j-4}}{(j-4)!} & \cdots & 1 \\
                                  \frac{x_{j}^{j}}{j!}&  \frac{x_{j}^{j-1}}{(j-1)!} &
                                \frac{x_{j}^{j-2}}{(j-2)!} &  \frac{x_{j}^{j-3}}{(j-3)!} & \cdots & \frac{x_j}{1} \\
                              \end{array}
                            \right), j\ge2.
\EQNY
\ET

}

\section{Examples}
In this section, \ZX{we present  several illustrating} examples.



\textbf{Example 1}. \LLP{Let $X(t)=B_H(t)+B_{1/2}(t^{2H}), t\in[0,T]$, with $H\in(0,1/2)$. Assume that \fD{the trend function $\gTt$ satisfies \eqref{eqgg} with some  constant $\mathcal{M}$ and some $d\ge1$.} We have
 \BQNY
\pb{\sup_{t\in[0,T]}\Bigl( X(t)-g(t) \Bigr)>u} \sim \frac{\mathcal{H}_{H}}{4^{\frac{1}{2H}}H}\Lambda_u^{1/H-2} \Psi\Bigl( \Lambda_u\Bigr),\ \  u\to\IF,
\EQNY
with
$$ \Lambda_u= \frac{u+ g(T)}{\sqrt{2}T^H}.$$
}
\COM{\sl Let $\{B_{{1}/{2}}^{(i)}(t^{2H_{i}}), t\in[0,T]\}, i\leq n,$ be independent scaled Brownian
 motions with parameters $H_{i}\in (0,1]$, respectively. 
Let $X_{i}(t)=B_{{1}/{2}}^{(i)}(t^{2H_{i}})$, $\lambda_{i}>0, \cH{i\le n}.$
 Then both the \aH{claims} of Theorem \ref{thm2} and \eE{Corollary 3.5(iii)} hold.
}
\COM{\textbf{Example 2}. {\sl Consider the fundamental martingales $\{M_{H_i}(t),t\in[0,T]\}, i\le n$, with $H_i\in(1/2,1)$, associated with independent  \LLP{fBm}'s  $\{B_{H_i}(t), t\in[0,T]\}, i\leq n,$ (see e.g.,
\aH{Norros} et al.\ (1999)) defined as
\BQNY
M_{H_i}(t)=\int_0^t\frac{s^{\frac{1}{2}-H_i}(t-s)^{\frac{1}{2}-H_i}}{2H_i \Gamma\left(\frac{3}{2}-H_i\right)\Gamma\left(\frac{1}{2}+H_i\right)}\, \aH{d B_{H_i}(s)},
\EQNY
\gT{with $\Gamma(\cdot)$ the Euler's Gamma function.}
Let $X_{i}(t)=M_{H_i}(t), i\le n$. Since $h^*(t)=\sum_{i=1}^n\frac{\lambda_i^2\Gamma\left(\frac{3}{2}-H_i\right)}{2H_i \Gamma\left({3}-{2}H_i\right)\Gamma\left(\frac{1}{2}+H_i\right)}t^{2-2{H_i}}$ is \LP{strictly }increasing\ZX{,} the
\aH{claim} of Theorem \ref{thm2} holds.
}}

\bigskip
The following time average Gaussian process was discussed in D\c{e}bicki and Tabi\'{s} (2011).

\textbf{Example 2}. { Let $\{B_{H_i}(t), t\in[0,T]\}, i\leq n,$ be independent \LP{\LLP{fBm's}} with Hurst parameters $H_{i}\in (0,1]$, $i\le n,$ \cl{satisfying} $H_{1}<H_{2}<\cdots<H_{n}$.
Set
\[
  X_{i}(t)=\left\{
 \begin{array}{cc}
  {\sqrt{2H_{i}+2}\frac{1}{t}\int_{0}^{t}B_{H_i}(s)ds},    & t>0,\\
  {0},     & t=0.
 \end{array}
  \right.
\]
Assume that \fD{the trend function $\gTt$ satisfies \eqref{eqgg} with some  constant $\mathcal{M}$ and some $d\ge1$.} It follows \cl{from Theorem \ref{thm1}} that
$$\mathbb{P}\left(\sup_{t\in[0,T]}\left( \ST- \gTt  \right)>u\right)
\sim \LLP{\Psi\left( \frac{u+ g(T)}{\cL{\sqrt{\sum_{i=1}^n \lambda_i^2 \fD{T^{2H_i}}}}} \right)}.$$
}




\textbf{Example 3}. \LP{Assume that $\LLP{U}(t)=\sum_{i=1}^{N(t)}Z_i, t\ge 0,$ is a compound Poisson process, with \LLP{i.i.d. claim} inter-arrival times  $\tau_i, i\inn$, being \cl{exponentially} distributed with parameter $\LLP{\mu}>0$, and}
  i.i.d claim sizes $Z_i, i\inn,$ having a Weibull distribution $F(y)=1-\exp(-y^\tau), \cH{y}\ge0,$ with shape parameter $\tau\in(0,1)$. Furthermore, let
 $X(t)=\sum_{i=1}^n \lambda_i \LLP{B_{1/2}(t^{2H_i})}$ with  $H_i\in(0,1], \cH{\lambda_i}>0,i\le n.$ In view of \nekorr{k4.2},
 we conclude that
\BQNY
\tilde{\psi}(u,c,T)\ \sim \ \LLP{\mu} T e^{-u^\tau},\quad \mathrm{as} \quad u \rightarrow \infty.
\EQNY
\LP{
\textbf{Example 4}. Consider a Gaussian perturbed $\alpha$-stable risk process.
Specifically, let $\{\LLP{U}(t), t\in[0,\IF)\}$ be an $\alpha$-stable L\'evy process with $\alpha\in(1,2)$,
i.e.  $\LLP{U}(t)\overset{d}=S_\alpha(t^{1/\alpha}, \beta,0)$, where
$S_\alpha(\sigma, \beta, d)$ denotes a stable random varible with index of stability $\alpha$,
scale parameter $\sigma$, skewness parameter $\beta$ and drift parameter $d$
(see e.g., Samorodnitsky and Taqqu (1994)).
Moreover, let $X(t)=\sum_{i=1}^n \lambda_i \cH{B}_{H_i}(t)$ with $B_{H_i}, i\le n,$ being independent \LLP{fBm}'s and $H_i\in(0,1], \cH{\lambda_i}>0,i\le n.$ It is known that  (\ref{eqpi}) is satisfied.
Consequently, it follows from \nekorr{k4.2} and the tail behavior of stable distribution (e.g., Samorodnitsky and Taqqu (1994)) that
\BQNY
\tilde{\psi}(u,c,T)\ \sim \ \pk{\LLP{U}(T)>u}\ \sim\ C_{\alpha,T^{1/\alpha}}\left(\frac{1+\beta}{2}\right)u^{-\alpha},\quad \mathrm{as} \quad u \rightarrow \infty,
\EQNY
where
\BQNY
C_{\alpha,T^{1/\alpha}}=\frac{T(\LPJ{1-\alpha})}{\Gamma(2-\alpha)\cos(\pi \alpha/2)}.
\EQNY
}

\COM{
\textbf{Example 6}. Consider a perturbed classical risk model with unit claims. Specifically, set
$$c(t)=ct, t\ge 0, \quad \delta_i=\delta, \quad Z_i\equiv1,i\inn,$$
 and define $X(t)=B_H(t)+B_{1/2}(t^{2H}), t\ge0$ with $H\in[\frac{1}{2},1)$. In addition, assume that $u\inn, u\ge2$, $f(t)=-\frac{c}{2}t$, $g(t)=\frac{c}{2}t$ and $\xi=\eta=1$. \gT{In view of \netheo{Th5.4} and  (\ref{eqb314}) we \ZX{have}}
\cH{the following bounds} for the \cH{finite-time} ruin probability
\begin{align*}
\left(1-\gT{{{\Upsilon}}_T}\left(u+1,\frac{3c}{2}\right)\right)&\left(1-\Psi\left(\frac{1+\frac{c T}{2}}{\sqrt{2T^{2H}}}\right)-\exp\left(\frac{-c T}{2T^{2H}}\right)\Psi\left(\frac{1-\frac{c T}{2}}{\sqrt{2T^{2H}}}\right)\right)\le \tilde{\psi}(u,c,T)\\
&\le 1-\gT{{{\Upsilon}}_T}\left(u-1,\frac{c}{2}\right)\left(1-\Psi\left(\frac{1+\frac{c T}{2}}{\sqrt{2T^{2H}}}\right)-\exp\left(\frac{-c T}{2T^{2H}}\right)\Psi\left(\frac{1-\frac{c T}{2}}{\sqrt{2T^{2H}}}\right)\right),
\end{align*}
with
\BQNY
\gT{{{\Upsilon}}_T}(u,r)=e^{-\delta T}\aH{\Biggl(}\sum_{m=0}^{[cT+u]}\frac{(\delta T)^m}{m!}+\sum_{j=u+1}^{[cT+u]}(-1)^{j+u}\left(\frac{\delta}{c}\right)^j(j-u)\left(\sum_{i=1}^{j-u}\frac{(-1)^{i+1}i^{j-1}}{(u+i)!(j-u-i)!}\right)\sum_{m=0}^{[cT+u]-j}\frac{(\delta T)^m}{m!}\aH{\Biggr)}.
\EQNY

}

\section{Proofs}
In this section we give detailed proofs of our previous results.
\ZX{Recall} that $X(t)=\sum_{i=1}^{n}\lambda_{i}X_{i}(t)$ is the \LLP{aggregate} \gT{ centered Gaussian process} with variance function
$\sigma_{X}^{2}(t):=\sum_{i=1}^{n}\cE{\lambda_{i}^2} \cE{\sigma_{i}^2(t)}$ and \aH{$\bar{X}(t):= X(t)/\sigma_X(t)$}.\\ 

\COM{
\prooftheo{thm2}
Since the processes $\XIT, i\le n,$ are independent, we have, for any $t\in [0,T]$, that
\BQNY
\mathbb{P}\left(\sup_{t\in[0,T]}\left(\sum_{i=1}^{n}\lambda_{i}X_{i}(t)- g(t)  \right)>u\right)& = &
\pk{\exists t \in [0,T]: \sum_{i=1}^{n}\lambda_{i}X_{i}(t)> u+g(t)}\\
&\ge & \pk{\sum_{i=1}^{n}\lambda_{i}X_{i}(t)> u+g(t)}\\
&=& \cT{\Psi\left(\frac{u+ g(t)  }{\sqrt{h^*(t)} }\right)},
\EQNY
\eE{and thus the first claim follows.}

We show the second claim by borrowing the idea of the proof of Theorem 2.4 of D\c{e}bicki and Sikora (2011).
\dT{Since the proofs for the three \gT{cases} are similar, we only give the one for   \LLP{sub-fBm} case.}
Define the Gaussian process
$$Z(t)=B_{1/2}(h^*(t)), \quad t\ge 0,$$
where \LP{$\{B_{1/2}(t), t\in[0,\infty\cE{)}\}$} is \eE{the} standard Brownian motion. \cE{For any $t\in [0,T]$,} we have 
$$\E{Z(t)}=0=\E{X(\cE{t})},\quad  \sigma_{Z}^{2}(t)=h^*(t)\LP{=\sum_{i=1}^{n}\lambda_{i}^{2}\LLP{(2-2^{2H_i-1})}t^{2H_{i}}}
.
$$
\dT{By the assumption that the process $ \XIT$ is a  \LLP{sub-fBm} with $H_i\in(0,1/2\gT{]}$, \LP{it is easy to see that}}
\LLP{\begin{eqnarray*}
\E{X(s)X(t)}
&=&\sum_{i=1}^{n}\lambda_{i}^{2}[\E{X_{i}(t)
X_{i}(s)}]
=\sum_{i=1}^{n}\lambda_{i}^{2}\left[s^{2H_{i}}+t^{2H_{i}}-\frac{1}{2}((t+s)^{2H_i}+\abs{t-s}^{2H_{i}})\right]\\
&\geq& \sum_{i=1}^{n}\lambda_{i}^{2}(2-2^{2H_i-1})\min(s^{2H_{i}},t^{2H_{i}})
= \E{Z(s)Z(t)}.
\end{eqnarray*}}
\dE{Hence, as in the aforementioned paper, applying Slepian's inequality (see e.g., Piterbarg (1996)),
 we obtain }
 \BQNY
\mathbb{P}\left(\sup_{t\in[0,T]}\left(X(t)- g(t)  \right)>u\right)& \le  &
1- \pk{\forall s \in [0,h^*(T)]: B_{1/2}(s)\le  u+g(h(s))}\\
&\le & 1- \pk{\forall s \in [0,h^*(T)]: B_{1/2}(s)\le  u+\overline{g_k}(h(s))},
\EQNY
with $\overline{g_k}$ being the polygon line with $k$ nodes $(t_0,g(h(t_0))) \ldot (t_k,g(h(t_k)))$ where $0=t_0 < t_1 < \cdots < t_k=h^*(T)$.
Note that the last bound follows by the assumption that $g(h(t))$ is a concave function, which implies that $u+ g(h(t))$ is also concave in $t$. In view of Janssen and Kunz (2004) (see also Hashorva (2005)), we have
\def\x{\vk{x}}
\BQNY
\lefteqn{\pk{\forall s \in [0,h^*(T)]: B_{1/2}(s)\le  u+\LLP{\overline{g_k}}(h(s))}}\\
& = & \Bigl[(2\pi)^{k} \prod_{i=0}^{k-1}(t_{i+1}-t_i) \Bigr] ^{-1/2}
   \int_{x_i < u+ g(h(t_i)), \forall i\le k}\exp(-\x^\top \Sigma^{-1} \x/2)\notag \\
&& \times\prod_{i=0}^{k-1}
\Bigl[1-\exp\bigl(-2( u+g(h(t_i)) - x_i)(u+g(h(t_{i+1}))-x_{i+1})/(t_{i+1}-t_i)\bigr)\Bigr]\,  d\x,
\EQNY
with $x_0=0, \x=(x_1 \ldot x_k)^{\cH{T}}$ and matrix $\Sigma^{-1}$ given by
\BQNY\label{eq:invSI}
  \hspace{-10mm}  \Sigma^{-1} = \left(
  \begin{array}{ccccc}
    \!\!\!  \frac{1}{t_1-t_0}+\frac{1}{t_2-t_1}  \!\!\! & -\frac{1}{t_2-t_1} &  0 & \cdots &  0 \\
     -\frac{1}{t_2-t_1} &  \!\!\!  \frac{1}{t_2-t_1}+\frac{1}{t_3-t_2}  \!\!\! &  -\frac{1}{t_3-t_2} & &  \vdots \\
     0 &  -\frac{1}{t_3-t_2} &  \!\!\!  \frac{1}{t_3-t_2}+\frac{1}{t_4-t_3} \!\!\! &  \ddots &  0 \\
     \vdots & &  \ddots &  \ddots &  -\frac{1}{t_k-t_{k-1}} \\
     0 & \cdots &  0 &  -\frac{1}{t_k-t_{k-1}} &  \!  \frac{1}{t_k-t_{k-1}}  \!\!\!
  \end{array} \right),
\EQNY
%
%
%
 %
%
 hence the proof follows easily. \QED\\

}

\prooftheo{thm1}
Define
$$ m_{u}(t):=\frac{u+ \gTt  }{\sigma_{X}(t)}\ \ \mbox{and}\quad \pi(u):=\mathbb{P}\left(\sup_{t\in[\eE{\delta},T]}\bar{X}(t)\frac{m_{u}(T)}{m_{u}(t)}>m_{u}(T)\right).
$$
\cE{For any $u\ge0$,} \aH{as in D\c{e}bicki and Sikora (2011),} 
we may further \ZX{write}
\BQNY
 \pi(u) \le \mathbb{P}\left(\sup_{t\in[0,T]}\left(X(t)- \gTt  \right)>u\right)
\leq \mathbb{P}\left(\sup_{t\in[0,\eE{\delta}]}\left(X(t)- \gTt  \right)>u\right)+\pi(u).
\EQNY
\COM{where
\begin{eqnarray*}
\pi(u)&=&\mathbb{P}\left(\sup_{t\in[\eE{\delta},T]}\left(X(t)- \gTt  \right)>u\right)\\
      &=&\mathbb{P}\left(\sup_{t\in[\eE{\delta},T]}\bar{X}(t)\frac{m_{u}(T)}{m_{u}(t)}>m_{u}(T)\right).
\end{eqnarray*}
}
\ZX{Obviously,}
$$1-\frac{m_{u}(T)}{m_{u}(t)}=\frac{\sigma_{X}(T)-\sigma_{X}(t)}
{\sigma_{X}(T)}+\frac{\sigma_{X}(t)[\gTt- \gTT]}{(u+ \gTt  )\sigma_{X}(T)}.$$
\fD{Further, in view of \eqref{eqgg}, $\delta$ can be suitably chosen such that
\BQNY
\big|g(T)-g(t)\big|\le Const(\sigma_{X}(T)-\sigma_{X}(t))
\EQNY
for all $t\in[\delta, T]$. Therefore, for any  $\varepsilon>0$, when $u$ is sufficiently large, we have, uniformly in $[\eE{\delta},T]$,}
\begin{eqnarray}
\label{eq3.1}
1-\cL{(1+ \varepsilon)}\frac{\sigma_{X}(T)-\sigma_{X}(t)}{\sigma_{X}(T)}
\leq\frac{m_{u}(T)}{m_{u}(t)}
\leq 1-(1-\varepsilon)\frac{\sigma_{X}(T)-\sigma_{X}(t)}{\sigma_{X}(T)}.
\end{eqnarray}
\fD{Consequently, it follows from (\ref{eq3.1}) that,} 
for $u$ sufficiently large,
$$\eE{ \pi_{+\ve}(u):=\mathbb{P}\left(\sup_{t\in[\eE{\delta},T]}Y_{+\ve}(t)>m_{u}(T)\right)}\le \pi(u)\leq \pi_{-\ve}(u):=\mathbb{P}\left(\sup_{t\in[\eE{\delta},T]}Y_{-\ve }(t)>m_{u}(T)\right),$$
%
where
$$\XE{ Y_{\pm \ve }\cE{(t)}:=\bar{X}(t)\left(1-(1 \pm \varepsilon)\frac{\sigma_{X}(T)-\sigma_{X}(t)}{\sigma_{X}(T)}\right),\ \ \fD{t\ge0}.}$$
\cE{Next, we \aH{analyse} $\pi_{-\ve }(u)$ \ZX{for fixed} $\eE{\ve}\in(0,1)$}, \XE{the asymptotics of $\pi_{+ \ve}(u)$ follows with the same arguments}. 
\fD{Obviously, the standard} deviation function  $\sigma_{Y_{\eE{-\ve}}}(t)$ attains its
unique maximum over $[\eE{\delta}, T]$ at $t=T$, with $\sigma_{Y_{\eE{-\ve}}}(T)=1$. 
 \LLP{Further, by Assumption \cE{A1} (recall  that $\tilde \sigma_i = \sigma_{i}(T)$)},
\begin{eqnarray*}
\sigma_{Y_{\eE{-\ve}}}(t)
      &=&1-(1-\eE{\ve})\frac{\tN}{\sum_{i=1}^{n}\lambda_{i}^{2} \tDHi }\LLP{(T-t)}^{\cT{\beta}}+o(\LLP{(T-t)}^{\cT{\beta}})
\end{eqnarray*}
as $t\uparrow T,$ with
$$\cT{\widetilde{N}=\lim_{t\rightarrow T}\sum_{i=1}^{n} \lambda_{i}^{2}\widetilde{\sigma_{i}}\tAi\LLP{(T-t)}^{(\beta_{i}-\beta)}}\LLP{\in(0,\IF)},$$
and
\COM{
\begin{eqnarray*}
\CO{\bar{Y}_{\eE{\ve}}(s),\bar{Y}_{\eE{\ve}}(t)}
      &=&\frac{1}{\sigma_{X}(s)}\frac{1}{\sigma_{X}(t)}\sum_{i=1}^{n}\lambda_{i}^{2}\sigma_{i}(s)\sigma_{i}(t)
      -\frac{\dT{\widetilde{G}}}{\sum_{i=1}^{n}\lambda_{i}^{2} \tDHi }|t-s|^{\AL}+o(|t-s|^{\AL})\\
\end{eqnarray*}
as $s,t\rightarrow T$, with
\dT{$$\widetilde{G}=\lim_{t,s\rightarrow T}\sum_{i=1}^{n}\lambda_{i}^{2}\tDi\tDHi|t-s|^{(\alpha_{i}-\AL)}\in(0,\infty).$$}
Further, as $s,t\rightarrow T$ we may write 
\begin{eqnarray*}
\frac{\sum_{i=1}^{n}\lambda_{i}^{2}\sigma_{i}(s)\sigma_{i}(t)}{\sigma_{X}(s)\sigma_{X}(t)}
&=&1-\frac{ \tMT     }
{2(\sum_{i=1}^{n}\lambda_{i}^{2} \tDHi )^{2}}|t-s|^{2\BL }+o(|t-s|^{2\BL })
\end{eqnarray*}
with
$$ \tMT     =\lim_{\eE{s, t\rightarrow }T}\sum_{i=2}^{n}\sum_{j=1}^{i-1}\lambda_{i}^{2}\lambda_{j}^{2}( \tAj \cT{\widetilde{\sigma_{i}}}- \tAi  \cT{\widetilde{\sigma_{j}}}
|t-s|^{\beta_{i}-\beta_{j}})^{2}|t-s|^{2(\beta_{j}-\BL )}\in(0,\infty).$$

\cE{ Utilising the fact that both  $\AL$ and $\BL$ are positive,  we obtain that}

\[
  1-\CO{\bar{Y}_{-\eE{\ve}}(s),\bar{Y}_{-\eE{\ve}}(t)}=\left\{
 \begin{array}{cc}
  {K_{1}|t-s|^{\AL}+o(|t-s|^{\AL})} ,    & \AL< 2\BL ,\\
  {(K_{1}+K_{2})|t-s|^{\AL}+o(|t-s|^{\AL})} ,    & \AL= 2\BL ,\\
  {K_{2}|t-s|^{2\BL }+o(|t-s|^{2\BL })},    & \AL> 2\BL
 \end{array}\right.
\]
}
$$\LLP{
  1-\CO{\bar{Y}_{-\eE{\ve}}(s),\bar{Y}_{-\eE{\ve}}(t)}=
\frac{\dT{\widetilde{G}}}{\sum_{i=1}^{n}\lambda_{i}^{2} \tDHi } |t-s|^{\AL}+o(|t-s|^{\AL})
}$$
as $\LPJi{\min(s,t)}\rightarrow T$,
with
\COM{
$$K_{1}=\frac{\dT{\widetilde{G}}}{\sum_{i=1}^{n}\lambda_{i}^{2} \tDHi }, \quad K_{2}=\frac{ \tMT     }
{2(\sum_{i=1}^{n}\lambda_{i}^{2} \tDHi )^{2}},$$
with constants 
}
\dT{$$\widetilde{G}=\lim_{t,s\rightarrow T}\sum_{i=1}^{n}\lambda_{i}^{2}\tDi\tDHi|t-s|^{(\alpha_{i}-\AL)}\in(0,\infty).$$}

\LLP{Moreover, i}n view of Assumption \cE{A2}, \LP{we have,} \aH{for $s,t\in[\eE{\delta},T]$ and some $\mathbb{C}>0$}, 
\begin{eqnarray*}
\E{(Y_{-\eE{\ve}}(t)-Y_{-\eE{\ve}}(s))^{2}}
&=& \E{ \left(\eE{\ve}(\bar{X}(t)-\bar{X}(s))+\frac{1-\eE{\ve}}{\sigma_{X}(T)}(X(t)-X(s))\right)^{2}}\\
&\leq & 2\eE{\ve^2}\E{(\bar{X}(t)-\bar{X}(s))^{2}}+\frac{2(1-\eE{\ve})^{2}}{\sigma_{X}^{2}(T)}\E{(X(t)-X(s))^{2}}\\
&\leq & \left(\frac{2\eE{\ve}^{2}}{\sigma_{X}^{2}(\eE{\delta})}+\frac{2(1-\eE{\ve})^{2}}{\sigma_{X}^{2}(T)}\right)
\E{(X(t)-X(s))^{2}}\\
&\leq & \LLP{\mathbb{C}}|s-t|^{\min_{1\leq i\leq n}\gamma_{i}}.
\end{eqnarray*}
\cE{Therefore, the} \aH{Gaussian} process $\cH{\{\eE{Y}_{\eE{-\ve}}(t), t\in [0, T]\}}$ satisfies \LLP{the
conditions of Theorem 8.2 of Piterbarg (1996)}  with $$A=\left((1-\eE{\ve})\frac{\tN}{\sum_{i=1}^{n}\lambda_{i}^{2} \tDHi }\right)^{1/\beta},\quad\!\! \LLP{C=\left(\frac{\dT{\widetilde{G}}}{\sum_{i=1}^{n}\lambda_{i}^{2} \tDHi }\right)^{1/\alpha},}$$
\COM{
D=\left\{
\begin{array}{cc}
{K_{1}},& \AL< 2\BL ,\\
{K_{1}+K_{2}},& \AL= 2\BL ,\\
{K_{2}},& \AL> 2\BL,
\end{array}\right.
\dT{\widetilde{\alpha}}=\left\{
\begin{array}{cc}
\eE{{\AL}},& \AL\leq 2\BL ,\\
{2\BL },& \AL>2\BL,
\end{array}\right.
and $\dT{\widetilde{\beta}}=\beta$.
}
 \aH{\LLP{and thus,} as $u\to \infty$,}\\ 

\[
  \pi_{-\eE{\ve}}(u)\sim\left\{
 \begin{array}{cc}
 \mathcal{H}_{{\AL}/{2}}\eE{\beta^{-1}}\Gamma(1/\beta)A^{-1} C  \tHu \tPSu ,    & \AL< \BL ,\\
  \eE{\mathcal{P}_{\alpha}^{(1-\varepsilon)\widetilde{N}/\widetilde{G}}}\tPSu ,    & \AL= \BL ,\\
  \LP{\tPSu},    & \AL> \BL,
 \end{array}\right.
 \ \ \ \mathrm{with}\ \ \ \LLP{\tPSu:=\Psi\left( \frac{u+ g(T)}{\cL{\sqrt{\sum_{i=1}^n \lambda_i^2 \tSIi}}} \right).}
\]
Consequently, letting $\ve\to 0,$ 
\[
  \pi(u)\sim\left\{
 \begin{array}{cc}
 \mathcal{H}_{{\AL}/{2}}\eE{\beta^{-1}}\Gamma(1/\beta) C \left(\frac{\tN}{\sum_{i=1}^{n}\lambda_{i}^{2} \tDHi }\right)^{-1/\beta} \tHu \tPSu ,    & \AL< \BL ,\\
  \eE{\mathcal{P}_{\alpha}^{\widetilde{N}/\widetilde{G}}}\tPSu ,    & \AL= \BL ,\\
  \LP{\tPSu},    & \AL> \BL,
 \end{array}\right.
\]
as $u\to \infty$. \LLP{Finally, using Borell-TIS inequality (e.g., Adler and Taylor (2007)) we conclude, as $u\to \infty$,
\BQNY
\mathbb{P}\left(\sup_{t\in[0,\eE{\delta}]}\left(X(t)- \gTt  \right)>u\right)&\le&\mathbb{P}\left(\sup_{t\in[0,\eE{\delta}]}X(t)>u+\inf_{t\in[0,\delta]}g(t)\right)\\
&\le& \exp\left(\frac{-\left(u+\inf_{t\in[0,\delta]}g(t)-\mathbb{E}\left(\sup_{t\in[0,\eE{\delta}]}X(t)\right)\right)^2}{2\sigma_\delta^2}\right)=o(\pi(u)),
\EQNY
\fD{since} $\sigma_\delta^2:=\sup_{t\in[0,\delta]}(\sum_{i=1}^{n}\lambda_{i}^{2}\sigma_{i}^{2}(t))<\sum_{i=1}^{n}\lambda_{i}^{2}\widetilde{\sigma_{i}}^{2}$. The proof is complete.
 } \QED\\

\cl{The next lemma is crucial for \ZX{the proof of \nekorr{korr:1}. Details of its proof are omitted here} since there are only some algebra calculations involved.}

\BEL \label{lem:corr2}
Under the conditions of \nekorr{korr:1}, for any $i\le n$ and  $T>0$, \LLP{we have, as $s,t \to  T$}

\eE{(i) if $\XIT$ is a  \aH{bi-\LLP{fBm}} , then}
\BQNY
\sigma_{i}(t)&=&T^{K_iH_{i}}-K_i H_{i}T^{K_iH_{i}-1}\LLP{(T-t)}+o(\LLP{(T-t)}),\\
1-\CO{\bar{X}_{i}(t),\bar{X}_{i}(s)}& =& \frac{1}{2^{K_i} T^{2K_i H_{i}} }|t-s|^{2K_iH_{i}}+o(|t-s|^{2K_iH_{i}});
\EQNY

\eE{(ii) if $\XIT$ is a sub-\LLP{fBm}, then
\BQNY
\sigma_{i}(t)&=&\sqrt{2-{2^{2H_i-1}}}T^{H_{i}}-\sqrt{2-{2^{2H_i-1}}}H_{i}T^{H_{i}-1}\LLP{(T-t)}+o(\LLP{(T-t)}),\\
1-\CO{\bar{X}_{i}(t),\bar{X}_{i}(s)}& =& \frac{1}{2(2-2^{2H_i-1}) T^{2H_{i}} }|t-s|^{2H_{i}}+o(|t-s|^{2H_{i}});
\EQNY}
\COM{Furthermore,
if $\XIT$ is \cL{a} \LLP{fBm},
\BQNY
1-\CO{\bar{X}_{i}(t),\bar{X}_{i}(s)}& =& \frac{1}{2 T^{2H_{i}} }|t-s|^{2H_{i}}+o(|t-s|^{2H_{i}});
\EQNY
}
\COM{
\dT{(iii) if $\XIT$ is a  $H$-ss Gaussian martingale, then
\BQNY
\sigma_{i}(t)&=&\LP{\sqrt{\E{X_i^2(1)}}}T^{H_{i}}- \LP{\sqrt{\LLP{\E{X_i^2(1)}}}} H_{i}T^{H_{i}-1}\LLP{(T-t)}+o(\LLP{(T-t)}),\\
1-\CO{\bar{X}_{i}(t),\bar{X}_{i}(s)}& =& \frac{H_{i}}{T}|t-s|+o(|t-s|).
\EQNY}
}
Additionally, the \dE{process $\XIT$} satisfies the condition of Assumption \LLP{A2} for some \ZX{ positive} $\delta, \mathbb{C}$, and $\gamma_i=2K_i H_i$ and $H_i/2$ for  \aH{bi-\LLP{fBm}} and sub-\LLP{fBm}, respectively.
\EEL

\proofkorr{korr:1} The claim follows \LP{from} \netheo{thm1} and \nelem{lem:corr2},
where $\BL:=\beta_{1} =\beta_{2}=\cdots=\beta_{n}=1$, $ \tAi  =K_i H_{i}T^{K_i H_{i}-1}$, $\alpha_{i}=2K_iH_{i}$ and
$\LPJi{D_{i}}=\frac{1}{2^{K_i} T^{2K_iH_{i}} }$ for the  \aH{bi-\LLP{fBm}}; $\BL :=\beta_{1}=\beta_{2}=\cdots=\beta_{n}=1$, $ \tAi  =\sqrt{2-2^{2H_{i}-1}}H_{i}T^{H_{i}-1}$, $\alpha_{i}=2H_{i}$ and
$\LPJi{D_{i}}=\frac{1}{2(2-2^{2H_{i}-1})T^{2H_{i}}}$ for \cH{the} sub-\LLP{fBm}. 
\QED\\

\COM{
(iii) for each $T>0$ and some constant $\delta>0$,
$$\E{ (X_{i}(t)-X_{i}(s))^{2}}\leq C|t-s|^{2H_{i}}$$
for all $s,t\in[\delta,T]$ and some constant $C>0$.

\eE{\textbf{Proof.} The conditions on the variance functions for the three cases are obvious. Furthermore,
\COM{
\cH{Since} as $t\rightarrow T$
\begin{eqnarray*}
\sigma_{i}(t)&=&t^{K_iH_{i}}=T^{K_iH_{i}}-(T^{K_iH_{i}}-t^{K_iH_{i}})\\
&=&T^{K_iH_{i}}-K_i H_{i}T^{K_iH_{i}-1}\LLP{(T-t)}+o(\LLP{(T-t)})
\end{eqnarray*}
for \cH{the}  \aH{bi-\LLP{fBm}} . Similar argument yields the corresponding ones for \cH{the} sub-\Fbm and  $H$-ss Gaussian martingale. }
 for \cH{the}  \aH{bi-\LLP{fBm}} , we have
\begin{eqnarray*}
1-\CO{\bar{X}_{i}(t),\bar{X}_{i}(s)}
&=&\frac{1}{2^{K_i}}\frac{1}{t^{K_iH_{i}}s^{K_iH_{i}}}
\left[|s-t|^{2K_iH_{i}}-(s^{2H_{i}}+t^{2H_{i}})^{K_i}+2^{K_i}s^{K_iH_{i}}t^{K_iH_{i}}\right]\\
&=&\frac{1}{2^{K_i}T^{2K_iH_{i}}}|t-s|^{2K_iH_{i}}+o(|t-s|^{2K_iH_{i}}),
\end{eqnarray*}
as $t,s\rightarrow T$.
\cL{Similarly, we have, for \cH{the} sub-\Fbm,
\begin{eqnarray*}
1-\CO{\bar{X}_{i}(t),\bar{X}_{i}(s)}
=\frac{1}{2(2-2^{2H_i-1})T^{2H_{i}}}|t-s|^{2H_{i}}+o(|t-s|^{2H_{i}}),
\end{eqnarray*}
}
\dT{and, for  $H$-ss Gaussian martingale,
\begin{eqnarray*}
1-\CO{\bar{X}_{i}(t),\bar{X}_{i}(s)}
&=&\frac{1}{2}\frac{1}{t^{H_{i}}s^{H_{i}}}
\left[|t^{2H_{i}}-s^{2H_{i}}|-(t^{H_{i}}-s^{H_{i}})^{2}\right]\\
&=&\frac{1}{T}H_{i}|t-s|+o(|t-s|),
\end{eqnarray*}
as $t,s\rightarrow T$. The last claim follows easily, since
$$\E{(X_{i}(t)-X_{i}(s))^{2}}\le2^{1-K_i}|t-s|^{2K_iH_{i}}
$$
for  \aH{bi-\LLP{fBm}}  (cf.\ Proposition 3.1 of Houdr\'{e} and Villa (2003)),
\cL{$$\E{(X_{i}(t)-X_{i}(s))^{2}}\le K_{\epsilon,T}|t-s|^{H_{i}-\epsilon}$$
with $K_{\epsilon,T}$ some constant and $\epsilon\in(0,H_i)$ for \cH{the} sub- \LLP{fBm} (see (2.6) of Bojdecki et al.\ (2004)),}
 and
\begin{eqnarray*}
\E{(X_{i}(t)-X_{i}(s))^{2}}&= & |t^{2H_{i}}-s^{2H_{i}}| \leq  2H_{i}\cL{\max(T^{2H_{i}-1},\delta^{2H_{i}-1})}|t-s|
\end{eqnarray*}
for the  $H$-ss Gaussian martingale, and thus the proof is complete. }}\QED

}


\COM{
\prooftheo{korr:1} \eE{The claim follows applying \netheo{thm1} and \nelem{lem:corr2}
where $\BL:=\beta_{1} =\beta_{2}=\cdots=\beta_{n}=1$, $ \tAi  =K_i H_{i}T^{K_i H_{i}-1}$, $\alpha_{i}=2K_iH_{i}$ and
$\widetilde{D_{i}}=\frac{1}{2^{K_i} T^{2K_iH_{i}} }$ for  \aH{bi-\LLP{fBm}} , $\BL :=\beta_{1}=\beta_{2}=\cdots=\beta_{n}=1$, $ \tAi  =\sqrt{2-2^{2H_{i}-1}}H_{i}T^{H_{i}-1}$, $\alpha_{i}=2H_{i}$ and
$\widetilde{D_{i}}=\frac{1}{2(2-2^{2H_{i}-1})T^{2H_{i}}}$ for \cH{the} sub-\LLP{fBm}, and
$\BL:=\beta_{1} =\beta_{2}=\cdots=\beta_{n}=1$, $ \tAi  =H_{i}T^{H_{i}-1}$, $\AL:=\alpha_1=\alpha_{2} =\cdots =\alpha_{n}=1$ and
$\widetilde{D_{i}}=\frac{H_{i}}{T}$
for  $H$-ss Gaussian martingale.} \QED
}


\aH{\prooftheo{Th4.1} } 
We first give the proof of the second tail equivalence of (\ref{eqTh4.1_2}).
It is easy to see that
\BQN \label{eqeq1}
\pk{\limitsupT{(\LLP{U}(t)-c(t))>u}}&\le&\pk{\limitsupT{\LLP{U}(t)}+\limitsupT{(-c(t))}>u}
\EQN
and thus
\BQN
\limitsup{u}\frac{\pk{\limitsupT{(\LLP{U}(t)-c(t))>u}}}{1-{F_1}(u)}&\le& 1.\label{eq_3}
\EQN
\ee{Further we can write}
\BQNY
\pk{\limitsupT{(\LLP{U}(t)-c(t))>u}}&\ge&\pk{\limitsupT{\LLP{U}(t)}-\limitsupT{c(t)}>u}\\
&\ge&\pk{\limitsupT{\LLP{U}(t)}-\limitsupT{c(t)}>u, \limitsupT{c(t)}\le d(u)}\\
&\ge&\pk{\limitsupT{\LLP{U}(t)}>u+d(u)}\pk{c(T)\le d(u)},
\EQNY
which together with (\ref{eq_2}) yields
\BQNY
\limitinf{u}\frac{\pk{\limitsupT{(\LLP{U}(t)-c(t))>u}}}{1-{F_1}(u)}&\ge& 1, \label{eq_4}
\EQNY
and thus \LPJi{the second tail equivalence of} (\ref{eqTh4.1_2}) is established. Since $F_1\in \heavyH$, it follows using similar arguments and Theorem 2.13 in  Foss et al. (2011) that
\BQNY
\limitsup{u}\frac{\pk{\limitsupT{(\LLP{U}(t)-c(t)+X(t))>u}}}{1-{F_1}(u)}&\le&  \limitinf{u}\frac{1-{(1-\psi(\cdot,c,T))*F_2}(u)}{\psi(u,c,T)}\frac{\psi(u,c,T)}{1-{F_1}(u)}\ =\ 1,
\EQNY
\LPJ{where $(1-\psi(\cdot,c,T))*F_2(u)$ denotes the convolution of distributions $1-\psi(u,c,T)$ and $F_2(u)$}.
Moreover
\BQNY
\pk{\limitsupT{(\LLP{U}(t)-c(t)+X(t))>u}}
&\ge&\pk{\limitsupT{(\LLP{U}(t)-c(t))}>u+d(u)}\pk{\limitsupT{(-X(t))}\le d(u)}\\
&\sim&\pk{\limitsupT{\LLP{U}(t)}>u+d(u)}\pk{\limitsupT{(-X(t))}\le d(u)}
\EQNY
as $u\rightarrow \IF,$ which together with (\ref{eq_2}) yields
\BQNY
\limitinf{u}\frac{\pk{\limitsupT{(\LLP{U}(t)-c(t)+X(t))>u}}}{1-{F_1}(u)}&\ge& 1. \label{eq_5}
\EQNY
Consequently, 
\BQNY
\pk{\limitsupT{(\LLP{U}(t)-c(t)+X(t))>u}}\ \sim\ {1-{F_1}(u)}
\EQNY

as $u\rightarrow\IF$, \ee{and thus the claim follows.} \QED\\

\proofkorr{k4.2} 
\ee{By} Theorem 4.1 of Albin and Sund\'{e}n (2009)
$\LLP{U}(T)\in \longT$ and $\limitsupT{\LLP{U}(t)} \in \longT$ are equivalent.
\COM{ This implies
that there exists a function $d(u)$ such that $d(u)\rightarrow \IF$  
as $u\rightarrow\IF$, and
\BQNY
1-{F_1}(u\pm d(u))\ \sim\  1-{F_1}(u)\quad \mathrm{as}\ u\rightarrow\IF.
\EQNY}
\LLP{Consequently, the claim follows applying} Theorem \ref{thm1} and Theorem \ref{Th4.1}. \QED\\

\COM{

\aH{\prooftheo{Th5.4}}
\cL{By the assumption that $\{X(t), t\ge0\}$ and $\{S(t),t\ge0\}$ are independent, we see that
\BQNY
\pk{\US{t\in[0,T]}\tilde{S}(t)>u}&=&1-\pk{\US{t\in[0,T]}\tilde{S}(t)\le u}\\
&\le& \cl{1-\pk{\US{t\in[0,T]}(\LLP{U}(t)-(c(t)-g(t)))+\US{t\in[0,T]}(X(t)-g(t))\le u}}\\
&\le& 1-\pk{\US{t\in[0,T]}(\LLP{U}(t)-(c(t)-g(t)))\le u-\xi}\pk{\US{t\in[0,T]}(X(t)-g(t))\le \xi},
\EQNY
and
\BQNY
\pk{\US{t\in[0,T]}\tilde{S}(t)>u}&\ge& \cl{\pk{\US{t\in[0,T]}(\LLP{U}(t)-(c(t)-f(t)))-\US{t\in[0,T]}(-X(t)+f(t))> u}}\\
&\ge& \pk{\US{t\in[0,T]}(\LLP{U}(t)-(c(t)-f(t)))> u+\eta}\pk{\US{t\in[0,T]}(-X(t)+f(t))\le \eta}
\EQNY
for any $g\in \mathcal{G}_T$, $f\in \mathfrak{G}_T$, $\xi\in(0,u)$ and $\eta\in(0,\infty)$. The rest of the proof follows \cH{from} Theorem 3.4 and Eq. (3) of Ignatov and Kaish\cH{e}v (2000). \QED
}
}

\textbf{Acknowledgments.}  \clp{We would like to thank the referee and the editor for their comments and suggestions. 
 The authors kindly acknowledge partial support by the project RARE -318984 (a Marie Curie IRSES FP7 Fellowship) 
 and SNSF Grant 200021-140633/1. 
  K. D\c ebicki has been also supported by NCN Grant No 2013/09/B/ST1/01778 (2014-2016), and 
 Z. Tan from the NSF  of China (No. 11326175) and NSF of Zhejiang Province of China (No. Q14A010038).}

\end{document}